\documentclass[12pt,twoside]{amsart}

\usepackage[matrix,arrow]{xy}

\author{Florin Ambro} 
\address{Department of Mathematics\\
The Johns Hopkins University\\
3400 N. Charles, Baltimore MD 21218}
\email{ambro@chow.mat.jhu.edu}

\newcommand{\Q}{{\mathbb Q}}
\newcommand{\Z}{{\mathbb Z}}

\newcommand{\R}{{\mathbb R}}

\newcommand{\cH}{{\mathcal H}}

\newcommand{\relint}{\operatorname{relint}}
\newcommand{\Reg}{\operatorname{Reg}}
\newcommand{\Supp}{\operatorname{Supp}}

\newcommand{\mult}{\operatorname{mult} }
\newcommand{\codim}{\operatorname{codim} }

\newcommand{\calMld}{\operatorname{\mathcal Mld}}

\theoremstyle{plain}
\newtheorem{thm}{Theorem}[section]
\newtheorem{mthm}{Main Theorem}

\newtheorem{lem}[thm]{Lemma}
\newtheorem{cor}[thm]{Corollary}

\newtheorem{prop}[thm]{Proposition}

\newtheorem{conj}[thm]{Conjecture}
\newtheorem{hpth}[thm]{Hypothesis}

\theoremstyle{definition}
\newtheorem{defn}[thm]{Definition}

\newtheorem{notation}[thm]{Notation}

\newtheorem{exmp}[thm]{Example}

\newtheorem{rem}[thm]{Remark}

\newtheorem{ack}{Acknowledgments}   

\theoremstyle{remark}

\newcommand{\orb}{\operatorname{orb}}

\thanks{This work was partially supported by NSF Grant 
DMS-9800807}

\begin{document}

\bibliographystyle{amsalpha+}

\title[minimal log discrepancies]
{On minimal log discrepancies}
\maketitle

\begin{abstract}
 An explanation to the boundness of minimal log discrepancies conjectured by 
 V.V Shokurov would be that the minimal log discrepancies of a variety in
 its closed points define a lower semi-continuous function. We check this
 lower semi-continuity behaviour for varieties of dimension at most $3$ and 
 for toric varieties of arbitrary dimension.
\end{abstract}


\setcounter{section}{-1}


\section{Introduction}


\footnotetext[1]{1991 Mathematics Subject Classification. Primary: 14B05,
Secondary: 14E30.}

The Logarithmic Minimal Model Program (LMMP for short) predicts that an
algebraic variety can be simplified by performing a finite sequence of 
surgery operations (extremal contractions and flips). Although singularities 
appear naturally in the process, there exists a class of mild singularities 
preserved by these operations. It is expected that varieties with {\em only 
log canonical singularites} form the largest class in which LMMP works.
\par
These mild singularities are controlled by {\em minimal log discrepancies}
(m.l.d.'s for short), invariants introduced by V.V. Shokurov \cite{problems}. 
For instance, the m.l.d. of a variety $X$ in a nonsingular (Grothendieck) 
point $\eta \in X$ is just the codimension of $X$ in $\eta$.
\par
Related to the existence and termination of flips is the A.C.C. 
Conjecture, proven in codimension two \cite{term, acc}, and for 
$\Gamma=\{0\}$ in the case of toric varieties \cite{Bor} (see Section 1 
for definitions and notations):

\begin{conj}\label{acc_conj} \cite{problems}
Let $(X,B)$ be a log variety, and let
$\Gamma \subset [0,1]$ be a subset satisfying the descending chain condition.
Then the set
$$
A(\Gamma,n):=\{a(\eta;B); \codim (\eta, X)=n, b_j \in \Gamma \ \forall j\}
$$
satisfies the ascending chain condition (a.c.c. for short).
\end{conj}

Note that $A(\Gamma,n)$ satisfies a.c.c. iff it is bounded from above and 
it has no accumulation points from below. The following conjecture, proven
up to codimension three \cite{C3f,Mrk,Kaw}, proposes a 
sharp upper bound:

\begin{conj}\label{max_mld} \cite{problems}
Let $(X,B)$ be a log variety and let $\eta \in X$ be a Grothendieck point. 
Then the following inequality holds:
$$a(\eta;B) \le  \codim \eta.$$
Moreover, $X$ is nonsingular in $\eta$ if
$a(\eta;B)>\codim \eta -1$. 
\end{conj}

Our main interest is in the first part of Conjecture~\ref{max_mld}. 
We formulated a stronger form in \cite{Am}, as a lower semi-continuity
behaviour of minimal log discrepancies:

\begin{hpth}\label{lower_semicontinuity} 
Let $(X,B)$ be a log variety, and consider the function
$$
a: X \to \{-\infty \} \cup \R,\ \ x \mapsto a(x;B) 
$$
defined on the closed points of the variety $X$. Then $a$ is lower 
semi-continuous , i.e. every closed 
point $x \in X$ has a neighborhood $x \in U \subseteq X$ such that 
$$
a(x;B)=\inf_{x' \in U} a(x';B).
$$
\end{hpth}

We should note here that lower semi-continuity does not hold if we allow the 
codimension of the points to jump. It turns out that lower semi-continuity 
is in fact equivalent to the following stronger form of the inequality 
proposed in Conjecture~\ref{max_mld}:

\begin{hpth}\label{precise} Let $(X,B)$ be a log variety, 
and let $\eta,\xi \in X$ be two Grothendieck points such that 
$\eta \in \bar{\xi}$. Then
$$
a(\eta;B) \le a(\xi;B)+\codim(\eta,\xi).
$$
\end{hpth}

Hypothesis~\ref{precise} has interesting inductive properties (cf. Section 2).
Our main result is the following:

\begin{mthm} 
\begin{enumerate} 
  \item Hypotheses~\ref{lower_semicontinuity} 
               and ~\ref{precise} are equivalent. 
  \item Hypothesis~\ref{precise} is valid if one of the following 
     extra assumptions is satisfied:
          \begin{itemize}
                \item[a)] $\codim \eta \le 3$, or
                 \item[b)] $X$ is a torus embedding and $B$ is invariant 
                       under the torus action.
           \end{itemize}
\end{enumerate}
\end{mthm}

\par
In section $1$ we review basic definitions and results. The equivalence 
of Hypotheses~\ref{lower_semicontinuity} and ~\ref{precise} is
proved in section $2$, as a formal consequence of the Finiteness 
Theorem~\ref{finite}. The latter states that the set of all minimal log 
discrepancies of a log pair $(X,B)$ form a finite set ${\calMld}(X,B)$, 
called the {\em mld-spectrum} of $(X,B)$.
Moreover, the fibers of the {\em mld map} $a: X \to {\calMld}(X,B)$, 
defined on the closed points of $X$, give a finite partition of $X$ into 
constructible sets. 
Section $3$ contains the proof of Hypothesis~\ref{precise} under
the extra assumption $\codim \eta \le 3$. It is based on LMMP in
dimension $3$ (cf. \cite{mori-3folds, 3-flips, 3-models}), and on 
known results on $3$-dimensional canonical and terminal points 
(cf. \cite{C3f,canmod,mori-term,Mrk, Kaw,acc}). The last section 
is a good illustration for all the above: 
Hypothesis~\ref{lower_semicontinuity} follows from  
explicit formulae for minimal log discrepancies.

\begin{ack} I am grateful to Professor Vyacheslav V. Shokurov for
useful discussions and criticism.
\end{ack}


\section{Prerequisites}


A {\em variety} is a reduced irreducible scheme of finite type over
a fixed field $k$, of characteristic $0$. An {\em extraction} 
is a proper birational contraction of normal varieties. We will use
Zariski's Main Theorem in the following form: if $\mu:\tilde{X} \to X$
is an extraction and $x \in X$ is a closed point such that 
$\dim \mu^{-1}(x)=0$, then $\mu$ is an isomorphism over a neighborhood
of $x$ (cf. \cite[Exercises II.3.22, III.11.2]{Ha}).
\par
We denote by $\eta_X$ the generic point of a variety $X$. A Grothendieck 
point $\eta \in X$ is called {\em proper} if $\eta \ne \eta_X$. A
neighborhood of $\eta$ in $X$ is an open subset $U \subseteq X$ such
that $\eta \in U$.

\begin{defn} A {\em log pair\/} $(X,B)$
  is a normal variety $X$ equipped with an $\R$-Weil divisor $B$ such that 
  $K+B$ is $\R$-Cartier. $B$ is called the {\em pseudo-boundary\/} of the 
  log pair. A {\em log variety\/} is a log pair $(X,B)$ such 
  that $B$ is an effective divisor. 
\end{defn}

\begin{defn}
  \begin{enumerate}
    \item A log pair $(X,B)$ has {\em log nonsingular support\/}
  if $X$ is nonsingular and $\Supp(B)$ is a divisor with normal
  crossings \cite[0-2-9]{KMM}.
     \item A {\em log resolution \/} of a log pair $(X,B)$
 is an extraction $\mu:\tilde{X} \to X$ such that $\tilde{X}$ is
 nonsingular and $\Supp(\mu^{-1}(B)) \cup Exc(\mu)$ is a 
 divisor with normal crossings.
  \end{enumerate}
\end{defn}

\begin{defn}
If $(X,B)$ is a log pair and $\mu:\tilde{X} \to X$ is an extraction,
there exists a unique divisor $B^{\tilde{X}}$ on $\tilde{X}$ such that
 \begin{itemize}
   \item[i)] $B^{\tilde{X}}=\mu^{-1}B$ on $\tilde{X} \backslash Exc(\mu)$,
   \item[ii)] $\mu^*(K+B)= K_{\tilde{X}}+B^{\tilde{X}}$.
 \end{itemize}
The divisor $B^{\tilde{X}}$, called the {\sl log codiscrepancy divisor} of 
$K+B$ on $\tilde{X}$, determines a log pair structure on $\tilde{X}$. 
\par
The induced log pair $(\tilde{X},B^{\tilde{X}})$ has log nonsingular support
if $\mu:\tilde{X} \to (X,B)$ is a log resolution.
In the sequel, when we say that $\mu:(\tilde{X},\tilde{B}) \to (X,B)$
is a log resolution, it is understood that $\tilde{B}=B^{\tilde{X}}$.
\end{defn}

\begin{defn}
Let $(X,B)$ be a log pair. Let 
$E \subset \tilde{X} \stackrel{\mu}{\to} X$ be a 
prime divisor on an {\em extraction} of $X$. 
The {\em log discrepancy} of $E$ with respect
to $K+B$ (or with respect to $(X,B)$), is defined as
$$
a_l(E;X,B)=1-e
$$
where $e$ is the coefficient of $E$ in the log codiscrepancy divisor 
$B^{\tilde{X}}$.
By definition, $a_l(E;X,B)=1$ if $E$ is not in the support of $B^{\tilde{X}}$.
The {\em center of $E$ on $X$} is $\mu(E)$, denoted by $c_X(E)$.
The log discrepancy $a_l(E;X,B)$ depends only on the discrete 
valuation defined by $E$ on $k(X)$, in particular independent
on the extraction $\tilde{X}$ where $E$ appears as a divisor. 
\newline
We will write $a(E;X,B)$ or $a(E;B)$, dropping 
the index $l$ and even the variety $X$ from the notation. 
However, $a(E;B)$ should not be confused with the standard notation
in the literature for the {\em discrepancy} of $K+B$ in $E$, which 
is equal to $-1+a_l(E;X,B)$.
\end{defn}

\begin{rem} In the above notation, the 
log discrepancies for prime divisors on $\tilde{X}$ are 
uniquely determined by the formula
$$
\mu^*(K_X+B)=K_{\tilde{X}}+\sum_{E \subset \tilde{X} } (1-a(E;X,B))E
$$
where the sum runs over all prime divisors of $\tilde{X}$. 
\end{rem}

\begin{defn} (V.V. Shokurov) The {\it minimal log discrepancy} of a log pair
$(X,B)$ in a proper Grothendieck point $\eta \in X$ is defined as
$$
a(\eta;X,B)=\inf_{c_X(E)=\eta} a(E;X,B),
$$
where the infimum is taken after all prime divisors on extractions of $X$
having $\eta$ as a center on $X$. We set by definition $a(\eta_X;X,B)=0$.
\end{defn}

\begin{defn} The log pair $(X,B)$ has only {\em log canonical 
 (Kawamata log terminal) singularities} if 
$a(\eta;B) \ge 0$ ($a(\eta;B)> 0$) for every proper point $\eta \in X$.
Also, $(X,B)$ is said to have only {\em canonical (terminal) singularities} 
if $a(\eta;B) \ge 1$ ($a(\eta;B)> 1$) for every point $\eta \in X$ of
codimenision at least $2$.
\end{defn}

\begin{rem}\cite[0-2-12]{KMM} One can read the singularity type on a 
resolution. Indeed, assume $\mu:(\tilde{X},\tilde{B}) \to (X,B)$ 
is a log resolution. Then $(X,B)$ has only log canonical singularities 
(Kawamata log terminal singularities) iff the same holds for 
$(\tilde{X},\tilde{B})$. Since $(\tilde{X},\tilde{B})$ has log nonsingular 
support, this is equivalent to the fact that all the coefficients of 
$\tilde{B}$ are at most $1$ (strictly less than $1$). 
\end{rem}

Since any extraction is an isomorphism up to codimension $1$, one can
easily compute minimal log discrepancies in these cases. Indeed,
if $\codim \eta =0$, then $a(\eta;B)=0$ by definition.
If $\codim \eta =1$, then $a(\eta;B)=1-b_\eta$, where $b_\eta$ is the 
coefficient of $B$ in $\bar{\eta}$ (which is zero if $\eta$ is not in the 
support of $B$).
\par
In Grothendieck points of codimension at least $2$, the minimal log 
discrepancy is either a non-negative real number,
or $-\infty$:

\begin{prop}\label{ml} \cite[17.1.1]{Ko1} Let $(X,B)$ be a log
pair and let $\eta \in X$ with $\codim \eta \ge 2$.
\begin{enumerate}
  \item If $(X,B)$ is not log canonical in any neighborhood of $\eta$,
then 
$$a(\eta;B)=-\infty.$$
  \item Assume that $(X,B)$ is log canonical in a neighborhood of $\eta$.
Let $(\tilde{X},\tilde{B})$ be a log resolution of $(X,B)$ such that 
$\mu^{-1}(\bar{\eta})$ is a divisor and 
$\mu^{-1}(\bar{\eta}) \cup \Supp(\tilde{B})=\sum_i E_i$ 
has normal crossings. Then 
$$
a(\eta;B)= \min_{c_X(E_i)=\eta } a(E_i;B) \in \R_{\ge 0}
$$
\end{enumerate}
\end{prop} 

\begin{lem}\label{ml_lemma} Under the same assumptions, the following hold:
\begin{itemize}
   \item[a)] $a(\eta;B)=-\infty$ if $a(\eta;B)<0$.
   \item[b)] $a(\eta;B)=-\infty$ if $\eta \in E$ is a proper point of a
             prime divisor $E$ with $a(\eta_E;B)<0$.
\end{itemize}
\end{lem}

\begin{proof} (of Proposition ~\ref{ml})
$i):$ By Lemma~\ref{ml_lemma}$.a)$, we just need to show that $(X,B)$ is 
log canonical in some neighborhood of $\eta$ if $a(\eta;B) \ge 0$. 
Suffices to show that $a(\xi;B) \ge 0$
for all $\eta \in \bar{\xi}$. Assume by contradiction that $a(\xi;B)< 0$.
Let $E$ be a prime divisor on an extraction $\mu:\tilde{X} \to X$ such that 
$a(E;B)<0$ and $c_X(E)=\xi$.
Since $\eta \in \bar{\xi}$ is a proper point, there exists a
proper point $\eta' \in E$ such that $c_X(\eta')=\eta$. From 
Lemma~\ref{ml_lemma}$.b)$,
$a(\eta';B)=-\infty$, hence $a(\eta;B)=-\infty$. Contradiction!
\par
$ii):$ Follows from Lemma~\ref{formula}.
\end{proof}

\begin{proof} (of Lemma~\ref{ml_lemma}) $a):$ 
Let $E$ be a prime divisor on an
extraction $\mu:\tilde{X} \to X$ such that $a(E;B)<0$ and $c_X(E)=\eta$.
Since the induced map $\mu|_E:E \to \bar{\eta}$ has generic fibers
of positive dimension, there exists a proper point $\eta' \in E$ such that
$c_X(\eta')=\eta$. Thus $a(\eta;B) \le a(\eta';\tilde{B})$, hence suffices
to check $b)$.
\par
$b):$ We may assume that $(X,B)$ has log nonsingular support.
Let $E_1$ be the exceptional divisor on the blow-up in $\eta$, 
and let $\eta_1$ be a component of $E \cap E_1$ dominating $\eta$. 
Inductively, let $E_{k+1}$ be the exceptional divisor on the blow-up in 
$\eta_k$, and let $\eta_{k+1}$ be a component of $E \cap E_{k+1}$ dominating 
$\eta_k$. A simple computation gives
$$
a(E_{k+1};B)=k \cdot a(E;B)+a(E_1;B), \ \ c_X(E_k)=\eta \ \ \forall k.
$$ 
In particular, $a(\eta;B)=-\infty$.

\end{proof}
 
\begin{notation} Assume $X$ is a nonsingular variety and $\cup_{i \in I} E_i$
is a divisor with normal crossings supporting the divisor 
$B=\sum_{i \in I} (1-a_i)E_i$.
\begin{itemize}
 \item For $J \in {\mathcal P}(I)$ denote $E_J=\cap_{j \in J}E_j$ 
        and $a_J=\sum_{j \in J}a_j$ (set $E_{\emptyset}=X$
        and $a_{\emptyset}=0$);
 \item For $\eta \in X$ set $I(\eta)=\{i \in I;\eta \in E_i\} 
          \in {\mathcal P}(I)$, and let $C_\eta$ be the generic point of the
       unique component of $E_{I(\eta)}$ containing $\eta$.
\end{itemize}
\end{notation}

\begin{lem}\label{formula} Assume that $(X,B)$ is a log nonsingular pair 
having only log canonical singularities in $\eta \in X$. Then
$$
a(\eta;B)=a_{I(\eta)} +\codim\eta -|I(\eta)|.
$$
In particular, $a(\eta;B)=a(C_\eta;B)+\codim(\eta,C_\eta)$.
\end{lem}

\begin{rem} In other words, if non-negative, minimal log discrepancies on 
log pairs with log nonsingular support are attained on the first blow-up. 
This is definitely false in general.
\end{rem}

\begin{proof} 
$Step \ 1:$ We first check that $\inf_{\eta \in E_{i_0}} a(\eta;B)= a_{i_0}$.
Indeed, $a(\eta;\sum_i E_i)\ge 0$ \cite[0-2-12]{KMM}, so $a(\eta;B) \ge
a_{i_0}$ if $\eta \in E_{i_0}$. The equality is attained on the generic point
of $E_{i_0}$.
\par
$Step \ 2:$ Assume $\eta$ is the generic point of a connected component of 
$E_{I(\eta)}$. Blowing up $X$ in $\eta$ we have again a log nonsingular
pair, and the new divisor $E$ has log discrepancy $a_{I(\eta)}$. 
From the previous case, we infer $a(\eta)= a_{I(\eta)}$. 
\par
$Step \ 3:$ Otherwise, shrinking $X$, we may assume that 
$I(\eta)=I$, and there exist divisors $\{B_j\}_{j \in J}$
such that $(X,B+\sum_{j \in J} B_j)$ has log nonsingular support
and $\bar{\eta}$ is a connected component of $\cap_{i \in I \cup J} B_i$.
Set $a_i=1$ for all $i \in J$. From Step 2,
$a(\eta;B)=a(\eta;\sum_{i\in I\cup J} (1-a_i)B_i)=
 a_I+|J| = a_I+\codim\eta -|I|$.
\end{proof}

\begin{exmp} Hypothesis~\ref{precise} is valid if we further assume that
$(X,B)$ has log nonsingular support.
\end{exmp}

\begin{proof} Indeed, let $\eta,\xi \in X$ with  $\eta \in\bar{\xi}$. 
There is nothing
to prove if $a(\eta;B)=-\infty$, so we may assume that $(X,B)$ has only
log canonical singularites in $\eta$. Then  
$
a(\eta;B)-(a(\xi;B)+\codim(\eta,\xi))=a_J-|J|\le 0
$
where $J=I(\eta) \setminus I(\xi)$.
\end{proof}

Minimal log discrepancies behave well with respect to products.

\begin{defn} If $(X,B_X)$ and $(Y,B_Y)$ are two log pairs, we denote
by $(X\times Y,B_{X\times Y})$ the {\em product log pair}, i.e. the usual 
product with canonical Weil divisor $K_{X\times Y}=K_X\times Y+ X\times K_Y$
and pseudoboundary $B_{X\times Y}=B_X\times Y+X\times B_Y$. Note that
$$
K_{X\times Y}+B_{X\times Y}= p_1^*(K+B_X)+p_2^*(K+B_Y)
$$
where $p_1$ and $p_2$ are the projections.
\end{defn}

One can easily check that $(X\times Y,B_{X\times Y})$ has log nonsingular 
support if so do $(X,B_X)$ and $(Y,B_Y)$. Moreover, If 
$\mu:(\tilde{X},\tilde{B})\to (X,B)$ and 
$\nu:(\tilde{Y},\tilde{D})\to (Y,D)$ are log resolutions then
$$
\varphi=\mu\times \nu:(\tilde{X}\times \tilde{Y},
   \tilde{B}\times \tilde{Y}+\tilde{X}\times \tilde{D}) \to
               (X\times Y,B\times Y+X\times D)
$$
is a log resolution. We will need the folowing lemma:

\begin{lem}\label{mld_product} 
   Assume $\eta$ and $\xi$ are points on the log pairs $(X,B_X)$
   and $(Y,B_Y)$, respectively. Then
            $$a(\eta\times \xi;B_{X\times Y})=a(\eta;B_X)+a(\xi;B_Y).
      $$
\end{lem}

\begin{proof} For good resolutions that compute minimal log discrepancies,
$$B_{\tilde{X}\times \tilde{Y}}=\sum_i (1-a_i)E_i\times \tilde{Y}+
                  \sum_j(1-b_j)\tilde{X}\times F_j$$
is a divisor with normal crossings, and 
      $$ \varphi^{-1}(\eta \times \xi)=\bigcup_{c_X(E_i)=
                \eta,c_Y(F_j)=\xi} 
         (E_i\times \tilde{Y}\cap \tilde{X}\times F_j).$$ 
For simplicity, we may assume $a(\eta;B_X), a(\xi;B_Y) \ge 0$
(the other cases are similar).
Therefore $a_i,b_j \ge 0$ near $\eta$ and $\xi$ respectively, and 
Lemma~\ref{formula} gives
\begin{equation*}
\begin{split}
a(\eta\times \xi;B_{X\times Y})
&= \min_{c_X(E_i)=\eta,c_Y(F_j)=\xi} 
                (a_i+b_j)
\\
&=\min_i a_i +\min_j b_j 
\\
&=a(\eta;B_X)+a(\xi;B_Y).
\\
\end{split}
\end{equation*}

\end{proof}


\section{The mld stratification}


\begin{defn} Let $(X,B)$ be a log pair. The set 
$$
\calMld(X,B):=\{a(\eta;B);\eta \in X\} \subset \{-\infty\} \cup \R
$$
is called the {\em mld-spectrum} of $(X,B)$. 
The partition of $X$ given by the fibers of the map
$$a:X \to \{-\infty\} \cup \R, \ x \mapsto a(x;B),$$
defined on the closed points of $X$, is called the 
{\em mld-stratification} of $(X,B)$. 
\end{defn}

\begin{thm}\label{finite} (Finiteness) The mld-spectrum $\calMld(X,B)$ of a 
log pair is a finite set, and the mld-stratification is constructible, i.e.
all the fibers of the map $a$ are constructible sets.
\end{thm}

\begin{prop}\label{key_fin} Assume $W \subset X$ is a closed irreducible 
subvariety and $(X,B)$ is a log pair with only log canonical singularities 
in $\eta_W$.
Then there exists an open subset $U$ of $X$ such that
$U \cap W\ne \emptyset$ and 
$$
a(x;B)=a(\eta_W;B)+\dim W
$$ 
for every closed point $x \in W\cap U$. 
\end{prop}

\begin{proof} (of Theorem~\ref{finite})
Suffices to prove that $a|_W$ takes a finite number 
of values and its fibers are constructible subsets, for every closed 
subset $W \subseteq X$. There is nothing to prove if $\dim W=0$, 
so let $\dim W>0$. Let $W_0$ be an irreducible components of $W$. From 
Lemma~\ref{ml_lemma}$.b)$ and Proposition~\ref{key_fin}, there exists an open 
subset $U_0 \subset X$ such that $U_0 \cap W_0 \ne \emptyset$, 
$a|_{U_0 \cap W_0}$ is constant, and $U_0$ does not intersect the other 
irreducible components of $W$. Thus
$$
W=(W\setminus U_0) \sqcup (W_0 \cap U_0),
$$
and we are done by Noetherian induction.
\end{proof}

\begin{proof} (of Proposition~\ref{key_fin}) Let 
$\mu:(\tilde{X},\tilde{B}) \to (X,B)$ be a log
resolution with a normal crossing divisor $\cup_{i \in I}E_i$ on
$\tilde{X}$ supporting $\tilde{B}=\sum_i (1-a_i)E_i$ and
the {\em divisor} $\mu^{-1}(W)$. Shrinking $X$ near $W$, we may assume
$$
\mu^{-1}(W)=\bigcup_{i \in I_W} E_i
$$ 
for some subset $I_W \subseteq I$, and $\mu(E_i)=W$ for every $i \in I_W$. 
We may assume that $(X,B)$ has only log canonical singularities, and 
$\dim W>0$. Note that 
$a(\eta_W;B)=\min_{i \in I_W} a_i$.
\par
Removing from $X$ all
components of $\mu(E_J) \ (\forall J \subseteq I)$ that do not contain $W$,
we may assume that $W \subseteq \mu(C)$, or $W\cap \mu(C)=\emptyset$
for every (connected) component $C$ of $E_J \ (\forall J \subseteq I)$.
We call {\em relevant} those components $C$ with $W=\mu(C)$. 
The following hold:
\begin{itemize}
   \item[a)] $I(\eta_C)\cap I_W \ne \emptyset$ and 
             $\dim C=\dim \tilde{X}-|I(\eta_C)|$ for every relevant $C$.
   \item[b)] If $\eta \in \mu^{-1}(W)$, then $C_\eta$, the unique component 
             of $E_{I(\eta)}$ containing $\eta$, is relevant.
\end{itemize}
Since the generic fibers of the morphisms $\mu|_C:C \to W$ have expected 
dimension, there exists an open subset 
$U \subset X$ such that $W \cap U \ne \emptyset$ 
and
$$
\codim (C \cap \mu^{-1}(x))=|I(\eta_C)|+\dim W
$$
for every relevant $C$ and for every closed point $x \in W \cap U$.
\par
Let $x\in W\cap U$ and $\eta \in \mu^{-1}(x)$. Then  
$\codim \eta \ge \codim (C_\eta \cap \mu^{-1}(x))= |I(\eta)|+\dim W$, hence 
$$a(\eta;\tilde{B})= a_{I(\eta)}+\codim \eta - |I(\eta)|  
\ge a(I(\eta))+\dim W.$$ 
But $I(\eta)\cap I_W \ne \emptyset$ and all $a_i$'s are non-negative 
numbers, hence $a(I(\eta))\ge a(\eta_W)$. Thus $a(\eta)\ge a(\eta_W)+\dim W$. 
Taking infimum after all
$\eta$'s as above, we obtain
$$
a(x;B) \ge a(\eta_W;B)+\dim W.
$$
Finally, let $k \in I_W$ be an index such that $a(\eta_W;B)=a_k$.
Let $\eta$ be the generic point of an irreducible component of
$E_k \cap \mu^{-1}(x)$ of maximal dimension. Since $E_k$ is relevant, 
$\codim \eta=\dim W+1$. Moreover, $I(\eta)=\{k\}$ since 
$\dim W+1=\codim \eta \ge \codim C_\eta =\dim W+|I(\eta)|$.
Therefore $a(\eta)=a_k+\codim \eta -1=a(\eta_W)+\dim W$, and the above 
inequality is in fact an equality.
\end{proof}

\begin{rem} Fix a log pair $(X,B)$ and consider the function
$$
a:X \to \{-\infty\} \cup \R,\ x \mapsto a(x;B)
$$
\begin{itemize}
   \item[a)] The fiber $a^{-1}(-\infty)$ is closed. It is the 
 union of all closed subvarieties $W$ of $X$ such that 
 $a(\eta_W;B)=-\infty$. Its complement is the biggest open subset 
 of $X$ on which $(X,B)$ has log canonical singularities.
   \item[b)] $a^{-1}(\{-\infty\} \cup \R_{\le 0})$ is the complement of
 the biggest open subset of $X$ on which $(X,B)$ has Kawamata log
 terminal singularities. It is denoted $Nklt(X,B)$, or $LCS(X,B)$.
   \item[c)] The fiber $a^{-1}(\dim X)$ contains the open dense subset
             $\Reg(X)\setminus \Supp(B)$, and the converse inclusion
   should hold if $B$ is effective, according to the second part of 
   Shokurov's Conjecture.  
\end{itemize}
\end{rem}

\begin{rem} Note that Lemma~\ref{mld_product} implies that
$$
\calMld (X\times Y, B_{X\times Y})=\calMld (X,B_X) +\calMld (Y,B_Y).
$$
\end{rem}

The first part of Conjecture~\ref{max_mld} can be reduced to 
Hypothesis~\ref{lower_semicontinuity}. Indeed, we may assume
$\eta=\{x\}$ is a closed point. The function $a(x;B)$ would 
jump downwards only in special points, and it is constant equal to $\dim X$ on
an open dense subset of $X$. Therefore $\sup_{x \in X} a(x;B) = \dim X$. 

\begin{lem}\label{equiv} The two hypotheses~\ref{lower_semicontinuity} and 
~\ref{precise} are equivalent.
\end{lem}

\begin{proof} Assume Hypothesis~\ref{precise} is valid, and let $x \in X$ be a 
closed point. Using Theorem~\ref{finite}, we may shrink $X$
such that $x \in \bar{C}$ for every
irreducible component $C$ of the fibers of the map $a$. 
For $x' \in X$, there exists a $C$ such that $x' \in C$. Since 
$x \in \bar{C}$, we infer that $a(x;B) \le a(\eta_C;B)+\dim \eta_C$. But
$a(\eta_C;B)+\dim \eta_C=a(x';B)$, so we are done.
\par
Assume Hypothesis~\ref{lower_semicontinuity} is valid. According to 
Proposition~\ref{key_fin}, we may assume that $\eta=\{x\}$ is a closed point 
and $x \in \bar{\xi}$. Let $U_x$ be a neighborhood of $x$ such that
$a(x;B) \le a(x';B)$ for all $x' \in U_x$. Then 
$U_x \cap \bar{\xi} \subset \bar{\xi}$ is an open dense subset.
From Proposition~\ref{key_fin}, there exists some $x' \in U_x \cap \bar{\xi}$ 
such 
that $a(x';B)=a(\xi;B)+\dim \xi$. Therefore $a(x;B) \le a(\xi;B)+\dim \xi$. 
\end{proof}

Hypothesis~\ref{precise} has very strong inductive properties.
Denote by $\cH_c$ the Hypothesis~\ref{precise} with the extra 
         assumption $\codim \eta =c$. 
         Fix $\eta \in X$ a point of codimension $c$, that we may assume
         to be a closed point $x$. 
         \begin{itemize}
              \item Suppose that $\cH_{c'}$ is valid for $c'<c$. Then $\cH_c$
         for $x$ is equivalent to the following weaker version: if
         $C$ is a curve passing through $x$, then $a(x;B) \le a(\eta_C;B)+1$
               \item Suppose $\cH_{c'}$ is valid for $c' \le c$ and the 
         characterization of nonsingularity from Conjecture~\ref{max_mld} is 
         valid for $c'<c$. Then $a(x;B)>\dim X -1$ implies that $x$ is
         an isolated singularity. Thus the new case in each dimension is 
         that of isolated singularities.
         \end{itemize}

\medskip
\section{Lower semi-continuity up to codimension $3$}


By Lemma~\ref{equiv}, suffices to check Hypothesis~\ref{precise}. 
We think of $\eta \in X$ as being fixed, and we shrink $X$ to neighborhoods
of $\eta$ without further notice. 
\par
We may assume $a(\eta;B)> 1$, otherwise
there is nothing to prove. Therefore $(X,B)$ has only log canonical 
singularities by Proposition~\ref{ml}. In particular, the coefficients
of $B$ are non-negative numbers less than or equal to $1$ (note that 
$(X,B)$ might not be Kawamata log terminal). 
Minimal log discrepancies are invariant to cutting with generic 
hyperplane sections, hence for our purposes we can always assume that 
some fixed Grothendieck point is in fact closed.

We will need the following results:

\begin{lem}\label{close} \cite{acc} Assume $(X,B)$ is a log variety, and 
$X$ is nonsingular in $\eta$. Then the following hold:
\begin{itemize}
  \item[i)] $a(\eta;B) \le \codim \eta$.
  \item[ii)] $a(\eta;B) \ge \codim \eta -1$ iff $\mult_\eta B \le 1$ and 
$a(\eta;B)=\codim \eta -\mult_\eta B$.
\end{itemize}
\end{lem}

\begin{prop}\label{dimtwo} \cite[3.1.2]{Al} Assume $\eta \in X$ is a singular 
point of codimension $2$ on the log variety $(X,B)$. Then $a(\eta;B) \le 1$. 
Moreover, equality holds iff $\eta \notin \Supp(B)$ and $X$ has a Du Val 
singularity in $\eta$. 
\end{prop}
 
\begin{lem} Hypothesis~\ref{precise} is valid if $\codim \eta \le 2$.
\end{lem}

\begin{proof} Since $a(\eta;B) \le 1$ if $\codim \eta  \le 1$, 
we may assume $\codim \eta =2$. By assumption, $a(\eta;B)>1$, hence
Proposition~\ref{dimtwo} implies that $X$ is nonsingular in $\eta$ and
$$
a(\eta;B)= 2- \mult_\eta B
$$ 
Therefore $a(\eta;B) \le 2=a(\eta_X;B)+ \codim \eta$. If $\codim \xi =1$, 
decompose $B=b\cdot \bar{\xi}+B'$, with $0 \le b \le 1$ and 
$\xi \notin \Supp B'$. Then $a(\eta;B) \le 2-b=a(\xi;B)+1$.
\end{proof}

\begin{thm}\label{3known} Assume $X$ is a $3$-fold and $K_X$ is 
$\Q$-Cartier (we take $B=0$).
The following hold for a singular closed point $x \in X$:
\begin{itemize}
   \item[i)]\cite[2.2]{C3f} If $(x,X)$ is a canonical singularity of 
            index $1$, then either $a(x)=1$, or $X$ has a cDV singularity
            at $x$, i.e. there exists a hyperplane section $H \subset X$
            having a Du Val singularity in $x$.
   \item[ii)]\cite{Mrk} $a(x)=2$ if $X$ has a cDV singularity in $x$.
   \item[iii)]\cite{Kaw} $a(x)=1+\frac{1}{r}$ if $(x,X)$ is a terminal
            singularity of index $r$.
\end{itemize}
\end{thm}

\begin{rem} See also \cite{Msk} for upper bounds of minimal 
            log discrepancies of certain hypersurface singularities.
\end{rem}

\begin{cor} Assume $\eta \in X$ is a point of codimension $3$ 
on the log variety $(X,B)$ such that $a(\eta;B) > 2$. Then $X$
is nonsingular in $\eta$.
\end{cor}

\begin{proof} We may assume $\dim X=3$ and $\eta=\{x\}$ is a closed point.
\par
$Step \ 1:$ By Lemma~\ref{why}, $a(\eta_C;B)>1$ for every curve passing
through $x$. From the codimension $2$ case, $(X,B)$ has only terminal
singularities. 
\par
$Step \ 2:$ $X$ has $\Q$-factorial singularities. Indeed, 
from LMMP we can find a $\Q$-factorialization 
$\mu:(\tilde{X},\tilde{B}) \to (X,B)$, where $(\tilde{X},\tilde{B})$ is a 
log variety again. If $\dim \mu^{-1}(x)>0$, there exists 
$\eta \in \mu^{-1}(x)$ with $\codim \eta \le 2$, hence 
$a(\eta;\tilde{B})\le 2$ from the codimension $< 3$ cases. Then
$a(x;B) \le 2$. Contradiction!
Otherwise, $\dim \mu^{-1}(x)=0$. Zariski's Main Theorem implies that 
$\mu$ is an isomorphism over a neighborhood of $x$, hence $X$ is 
$\Q$-factorial. 
\par
$Step \ 3:$ Assume by contradiction that $x$ is a singular point.
Then it must be an isolated terminal point. From 
Theorem~\ref{3known}, $a(x;B) \le a(x;0)=1+\frac{1}{r} \le 2$, where $r$ is 
the index of $K_X$ at $x$. Contradiction! 
\end{proof}

\begin{lem}\label{why}
Assume $x \in W \subset X$, and $\dim X=3$. Assume
that either $\codim W=1$ and $a(\eta_W;B) \le 0$, or $\codim W=2$ and 
$a(\eta_W;B) \le 1$. Then $$a(x;B) \le a(\eta_W;B)+\dim W.$$
\end{lem}

\begin{proof} We may assume $a(x;B) \ge 0$ and $a(\eta_W;B)\ge 0$. 
\par
$Step \ 1:$ Assume $\codim W=1$ and $a(\eta_W;B)=0$. By easy divisorial 
adjunction, $a(x;B) \le a(x;B_{W^\nu})$, where $B_{W^\nu}$ is the different 
of $K+B$ on the normalization $W^\nu$ of $W$. The log variety 
$(W^\nu,B_{W^\nu})$ has dimension $2$, so $a(x;B_{W^\nu}) \le 2$.
\par
$Step \ 2:$ Assume $\codim W=2$ and $0 \le a(\eta_W;B) \le 1$. 
From LMMP, there exists a crepant extraction 
$\mu:(\tilde{X},\tilde{B}) \to (X,B)$ such that $\tilde{B}$ is 
effective and there exists a prime divisor $E$ on $\tilde{X}$ with
$\mu(E)=W$ and $a(\eta_E;\tilde{B})=a(\eta_W;B)$. Let $\eta$ be the 
generic point of a curve in the fiber of $\mu|_E:E \to C$ over $x$.  
From the codimension $2$ case, $a(\eta;\tilde{B}) \le a(\eta_E;\tilde{B})+1$. 
But $a(x;B) \le a(\eta;\tilde{B})$, so we are done.
\end{proof}

\begin{prop}\label{argum} Hypothesis~\ref{precise} holds if $\codim \eta=3$.
\end{prop}

\begin{proof} We may assume that $\eta$ is a closed point $x$ on the 
$3$-fold $X$. 
\par
$Step \ 1:$ Assume $\bar{\xi}$ is a curve $C$ passing through $x$.
From Lemma~\ref{why}, we may assume that $a(\eta_C;B)>1$. Then we may
also assume $a(x;B)>2$, hence $X$ is nonsingular in both 
$x$ and $\eta_C$. By Lemma~\ref{close}, $a(x;B)=3-\mult_x B$ and 
$a(\eta_C;B)=2-\mult_C B$. Therefore
$$a(x;B)-(a(\eta_C;B)+1)=\mult_C B -\mult_x B \le 0.$$
\par
$Step \ 2:$ Assume $\bar{\xi}$ is a surface $S$ passing through $x$.
Let $x \in C \subset S$ be a curve. Then $a(\eta_C;B) \le a(\eta_S;B)+1$ 
from the codimension $2$ case. From the previous step we get
$a(x;B) \le a(\eta_C;B)+1$, thus $a(x;B) \le a(\eta_S;B)+2$.
\end{proof}

The following characterization of cDV singularities is part of the folklore, 
but we include here a proof for completeness.

\begin{cor} Assume $(X,B)$ is a log variety and 
$\eta \in X$ is a point of codimension $3$. Then $a(\eta;B) =2$ iff
exactly one of the following holds:
\begin{itemize}
   \item[i)] $\eta \notin \Supp(B)$ and $X$ has a cDV singularity 
             in $\eta$ (i.e. 
             a cDV singularity after cutting $\bar{\eta}$ with 
             $\codim \eta$ general hyperplanes).
   \item[ii)] $X$ is nonsingular in $\eta$ and $\mult_\eta B=1$.
\end{itemize}
\end{cor}

\begin{proof} The second part follows from Lemma~\ref{close}, so we just have 
to prove $i)$. We may assume $\dim X=3$ and $\eta=\{x\}$ is a closed singular 
point.
\par
$Step \ 1:$ $(X,B)$ has only canonical singularities. Indeed,  
$$a(\eta_C;B)\ge a(x;B)-1=1$$ for every curve passing through $x$. From the 
codimension 
$2$ case, $(X,B)$ has only canonical singularities on $X\setminus \{x\}$. But
$a(x;B)=2$, hence we are done. 
\par
$Step \ 2:$ Assume that $B$ is $\R$-Cartier. Then $K_X$ is $\Q$-Cartier,
and let $r$ be the index of $X$ at $x$. Since $2\le a(x;B) \le a(x;0) \le 2$,
we infer that $B=0$ near $x$ and $a(x):=a(x;0)=2$. 
We just have to prove that $r=1$, since then $X$ has only canonical Gorenstein 
singularities, and therefore $x \in X$ is a cDV point due to 
Theorem~\ref{3known}$.i)$.
Note that if $x \in X$ is a terminal point, then $a(x)=1+\frac{1}{r}$ by 
Kawamata, hence $r=1$.
\newline
$X$ admits a terminal crepant extraction by \cite[0.6]{canmod}, i.e. there 
exists an extraction 
$
\mu:\tilde{X} \to X
$
such that $\tilde{X}$ has only terminal singularities and 
$\mu^*K_X=K_{\tilde{X}}$. 
Note that $2=a(x) \le a(\tilde{x})$ for every closed point 
$\tilde{x} \in \mu^{-1}(x)$. Thus the terminal subcase implies
that $K_{\tilde{X}}$ is Cartier. In particular, $K_X$ is Cartier near
$x$, i.e. $r=1$.
\par
$Step \ 3:$ Assume that $B$ is not $\R$-Cartier at $x$. We have to show that 
this is impossible.
From LMMP we can find a small extraction $\mu:\tilde{X} \to X$
such that $\tilde{X}$ is $\Q$-factorial. Let $\tilde{B}$ be the proper
transform of $B$. In particular, $\mu^*(K+B)=K+\tilde{B}$. 
If $\eta \in \mu^{-1}(x)$, then
$$
2=a(x;B) \le a(\eta;\tilde{B}).
$$
We have $\dim \mu^{-1}(x)>0$. Otherwise, Zariski's Main Theorem would
imply that $\mu$ is an isomorphism over a neighborhood of $x$. Thus $B$ is 
$\R$-Cartier, contradicting our assumption.  
\par
Therefore $\mu^{-1}(x)$ is a connected union of curves and 
$\tilde{X}$ has only cDV isolated singularities in $\mu^{-1}(x)$ from Step 2.  
Moreover, $\tilde{B}$ intersects $\mu^{-1}(x)$ in a finite set of 
points. Otherwise, if some curve $C$ over $x$ is included in 
$\Supp(\tilde{B})$, then $a(\eta_C;B)<2$, a contradiction.
\par
We arrive at the final contradiction with the following argument, 
kindly suggested by V. V. Shokurov: $-K_{\tilde{X}}$ is $\mu$-nef, 
but not $\mu$-trivial, since $\tilde{B}$ intersects the fiber $\mu^{-1}(x)$.
However, $\tilde{X}$ admits no flipping contraction since
its {\em difficulty} \cite{nonv} is $0$. Contradiction!
\end{proof}


\section{Toric minimal log discrepancies}


We refer the reader to \cite{Fu} for definitions and basic notations of
toric geometry.
Let $X=T_Nemb(\Delta)$ be a toroidal embedding, and let $\{B_i\}_{i=1}^r$
be the $T_N$-invariant divisors of $X$, corresponding to the primitive
vectors $\{v_i\}_{i=1}^r$ on the $1$-dimensional faces $\Delta$. 
Note first that
$
K+\sum_i B_i \sim 0
$
and $(X,\sum_i B_i)$ is log canonical (cf. \cite[3.1]{Al_lcsp}).
\par
Let $B= \sum_i (1-a_i)B_i$ be an invariant $\R$-divisor such that 
$K+B$ is $\R$-Cartier. This is equivalent to the existence of some linear
form $\varphi \in M_\R$ such that $\varphi(v_i)=a_i$ for every $i$. 
Moreover, assume $0 \le a_i \le 1$ for every $i$, hence $(X,B)$ is a 
log variety with log canonical singularities.

\begin{rem}\cite{Bor} Let $T_N(\Delta') \to X$ be the birational extraction 
induced by a subdivision $\Delta' \subset \Delta$, and let
$E_v \subset T_N(\Delta')$ be the
invariant prime divisor corresponding to a primitive vector $v \in |\Delta|$.
Then 
$$
a(E_v;B)=\varphi(v).
$$
Since any toric variety can be resolved by a basic subdivision of the fan,
we obtain the following formula for minimal log discrepancies in  
orbits:
$$
a_\sigma :=a(\eta_{\orb(\sigma)};B)=\inf\{\varphi(v); v \in \relint(\sigma)\}
,\ \sigma \in \Delta.$$
Here, $\relint(\sigma)$ denotes the 
relative interior of $\sigma \subset \R\sigma$, and $\orb(\sigma)$ is the 
$T_N$-orbit corresponding to the cone $\sigma \in \Delta$. We dropped the 
primitiveness assumption on the vectors since $\varphi$ is non-negative on 
$|\Delta|$. Note that $a_{\{0\}}=0$.
\end{rem}

\begin{prop}\label{mld_toric} In the above notations, let 
$X=\bigsqcup_{\sigma \in \Delta} \orb(\sigma)$ be
the partition of $X$ into $T_N$-orbits.
\begin{itemize}
   \item[i)] Each strata in the mld-stratification is a union of orbits.
          In other words, $a(x;B)=a_\sigma+\codim(\sigma)$ for every cone
          $\sigma \in \Delta$ and every closed 
          point $x \in \orb(\sigma)$.
   \item[ii)] $a_\sigma+\codim(\sigma) \le a_\tau+\codim(\tau)$ for all 
               cones $\tau, \sigma \in \Delta$ such that $\tau$ is a face of
               $\sigma$ (i.e. $\orb(\sigma)$ is in the closure of 
               $\orb(\tau)$).
\end{itemize}
\end{prop}

\begin{rem} In particular, Hypothesis~\ref{lower_semicontinuity} 
is valid for toric varieties.
\end{rem}

\begin{proof}
\par
 $i):$ The equality holds for the generic closed point $x \in \orb(\sigma)$ 
       from Proposition~\ref{key_fin}. This extends to all the points in 
       $\orb(\sigma)$ since $T_N$ acts transitively on orbits and leaves 
       the boundary fixed.
\par
 $ii):$ Let $\tau$ be a proper face of $\sigma$ and let $a_\tau=\varphi(v)$
        for some $v \in \relint(\tau)$. 
        We can find primitive 
        vectors $v_{i_1},\ldots, v_{i_c}\ (c=\codim(\tau,\sigma))$
        on the $1$-dimensional faces of $\sigma$ such that 
        $$w=v+v_{i_1}+\cdots+ v_{i_c} \in \relint(\sigma).$$ Therefore
        $a_\sigma \le \varphi(w)=
       \varphi(v)+a_{i_1}+\cdots+a_{i_c} 
          \le a_\tau+ \codim(\tau,\sigma)$.  
\end{proof}

\begin{rem}
Assume $\tau \prec \sigma \in \Sigma$ and 
$a_\tau+\codim(\tau,\sigma)= a_\sigma$.
Let $\tau \prec \gamma \prec \sigma$. Then 
$a_\sigma \le a_\gamma+\codim(\gamma,\sigma) \le 
a_\tau+\codim(\tau, \sigma)$.
Therefore $a_\gamma+\codim(\gamma,\sigma)= a_\sigma$.
\end{rem}

The second part of Conjecture~\ref{max_mld} has the following interpretation
on toric varieties:

\begin{prop} Let $\sigma \subset N_\R$ be a strongly rational polyhedral 
cone generated by the primitive vectors $v_1,\ldots,v_r \in N$. Assume
$\varphi \in M_\R$ is a linear form such that $0 \le \varphi(v_i)\le 1$ 
for every $i$, and let
$$
\varphi_\sigma:=\inf\{\varphi(v);v \in \relint(\sigma)\} 
$$
If $\varphi_\sigma>\dim\sigma-1$ then $\sigma$ is a nonsingular cone.
\end{prop}

\begin{rem} According to Proposition~\ref{mld_toric}, under 
the above assumptions we have
$$\dim\sigma-1<\varphi_\sigma \le \dim\sigma.$$
Moreover, 
$\varphi_\sigma=\dim\sigma$ iff $\varphi(v_i)=1$ for every $i$. Indeed, 
the same equality must hold for any proper face of $\sigma$, in particular 
for the $1$-dimensional rays of $\sigma$, hence 
$\varphi(v_i)=\varphi_{\R_{\ge 0} \cdot v_i}=1$ for 
every $i$.
\end{rem}

\begin{proof} We use induction on $n=\dim\sigma$. If $n=1$, there
is nothing to prove, so let $n \ge 2$. 
By Proposition~\ref{mld_toric}$.ii)$, every proper
face $\tau \prec \sigma$ has the same property with respect to 
$\varphi|_{M_\tau \otimes \R}$. By induction, all proper faces of 
$\sigma$ are nonsingular cones.
\par
$Step \ 1:$ Assume $\sigma$ is a simplicial cone, i.e. $r=n$.
It is known that $\sigma$ is nonsingular iff 
$$P_\sigma=\{\sum_{i=1}^n t_i v_i \in \sigma \cap N; 0 \le t_i<1 \ \ 
\forall i\} =\{0\}.$$
Assume $P_\sigma \ne \{0\}$. Since all proper faces are non-singular cones, 
$P_\sigma \cap \partial(\sigma)=\{0\}$. Therefore there exists
$v=\sum_{i=1}^n t_i v_i \in P_\sigma \cap \relint(\sigma)$. Then $0<t_i<1$
for every $i$, hence 
$\bar{v}=\sum_{i=1}^n (1-t_i)v_i \in P_\sigma \cap \relint(\sigma)$. Therefore
$2\varphi_\sigma \le \varphi(v+\bar{v})=\sum_{i=1}^n \varphi(v_i) \le n$.
This implies $\varphi_\sigma \le \frac{n}{2}\le n-1$, a contradiction.
Therefore $P_\sigma=\{0\}$, hence $\sigma$ is a nonsingular cone.
\par
$Step \ 2:$ If $\tau\prec \sigma$ is a face of codimension $1$ and
            $v_i \notin \tau$, then $\tau+\R_{\ge 0} v_i$ is a 
            nonsingular cone of dimension $n$. Indeed,
 let $\sigma'=\tau+\R_{\ge 0} v_i \subseteq \sigma$. 
 By assumption, $\sigma'$ is a simplicial cone of dimension $\dim\sigma$. 
 This also implies that $\relint(\sigma') \subseteq \relint(\sigma)$, hence 
 $\sigma'$ has the same property with respect to 
 $\varphi|_{M_{\sigma'} \otimes \R}$.
 Therefore $\sigma'$ is nonsingular from Step 1.
\par
$Step \ 3:$ We may assume that $r=n+1$. Indeed, if $r=n$ we
   are done from Step 1. Otherwise, $r \ge n+1$, and we show that
   this leads to contradiction. Let $\sigma' \subseteq \sigma$ be a cone 
   of dimension $n$ generated by $n+1$ of the vectors $v_i$'s. 
   Then $\sigma'$ has the same property with respect to $\varphi|_{M_{\sigma'}
   \otimes \R}$, since $\relint(\sigma') \subseteq \relint(\sigma)$.
   Therefore suffices to show that the case $r=n+1$ is impossible.
 \par
$Step \ 4:$ Assume $\tau \prec \sigma$ is a face of codimension $1$ and
       $v_i, v_j \notin \tau$. Then $v_i\pm v_j \in \Z \cdot (\tau \cap N)$.
   Indeed, let $\{v_{k_1},\ldots, v_{k_{n-1}}\}$ be the generators of
   $\tau$, which also form a basis of the lattice $\Z\cdot (\tau \cap N)$.
   From Step 2, $\{v_i,v_{k_1},\ldots, v_{k_{n-1}}\}$ and
   $\{v_j,v_{k_1},\ldots, v_{k_{n-1}}\}$ are both basis for the lattice
   $\Z\cdot (\sigma \cap N)$. The transition matrix has determinant $\pm 1$,
   hence the statement.
 \par
$Step \ 5:$ Let $\sigma$ be generated by $\{v_1,\ldots,v_{n+1}\}$. By 
   Step 2, we may assume that $\{v_1,\ldots,v_n\}$ is a basis of the lattice
   $\Z\cdot (\sigma \cap N)$, hence
   $$v_{n+1}=\sum_{i=1}^n r_i v_i,\ \ r_i \in \Z.$$
   We show that $r_i \in \{-1,0,1\}$ for every $i$. 
   At least one $r_i$ is positive (negative). Assume $r_i>0$. Then
   $\{v_k;k \notin \{i,n+1 \}\}$ generates a codimension $1$ 
   face, hence $$v_i \pm v_{n+1} \equiv 0 
   \mod \sum_{k \notin \{i,n+1 \}}  \Z\cdot v_k.$$ 
   On the other hand, $v_{n+1} \equiv r_i v_i \mod 
   \sum_{k \notin \{i,n+1 \}} \Z\cdot v_k$, hence $r_i=\pm 1$. Therefore
   $r_i=1$.
\par 
   Assume $r_j<0$. If $r_i>0$, then 
         $\{v_k;k \notin \{i,j\}\}$ generates a codimension $1$ face, so 
   $$v_i \pm v_j \equiv rv_{n+1} \mod 
         \sum_{k \notin \{i,j,n+1 \}} \Z\cdot v_k.$$ 
   Since
   $v_{n+1}\equiv v_i+r_jv_j \mod \sum_{k \notin \{i,j,n+1 \}}\Z\cdot v_k $, 
   we deduce that $r=1$. Therefore $$v_j \equiv v_{n+1} \mod
     \sum_{k \notin \{j,n+1 \}} \Z\cdot v_k.$$ But 
   $v_{n+1}\equiv r_jv_j \mod \sum_{k \notin \{j,n+1 \}}\Z\cdot v_k$, thus
   $r_j=\pm 1$. Therefore $r_j=-1$.  
\par
$Step \ 6:$ Let 
        $v_{n+1}=v_1+\ldots +v_s-v_{s+1}-\ldots -v_k$, where 
       $s\ge 1$ and $s+1 \le k \le n$. 
      One can easily check that 
      $$v=v_1+\ldots+v_s+v_{k+1}+\ldots+v_n \in \relint(\sigma).$$
      Therefore $\varphi_\sigma \le \varphi(v)=\sum_{i=1}^s \varphi(v_i)+
      \sum_{i=k+1}^n \varphi(v_i) \le s+n-k \le n-1$. 
      Contradiction!
\end{proof}


\end{document}